\documentclass[10pt]{article}
\usepackage{geometry}                
\usepackage{amsmath, amsthm, amssymb}
\usepackage{graphicx, color}
\usepackage{pdfsync}
\newtheorem{lemma}{Lemma}
\newtheorem{corollary}[lemma]{Corollary}

\newtheorem{remark}[lemma]{Remark}

\newtheorem{theorem}[lemma]{Theorem}
\newtheorem{definition}[lemma]{Definition}
\usepackage{graphicx}
\usepackage{amssymb}
\usepackage{epstopdf}
\newcommand{\e}{\varepsilon}

\def\wto{\rightharpoonup}

\title{Minimizing movements along a sequence of functionals\\ and curves of maximal slope}
\author{Andrea Braides \\ \small Dipartimento di Matematica, Universit\`a di Roma Tor Vergata
\\ \small  via della ricerca scientifica 1, 00133 Roma, Italy\\ \\Maria Colombo\\
 \small ETH Institute for Theoretical Studies\\ \small 
Clausiusstrasse 47
8092 Z\"urich, Switzerland\\
\\ Massimo Gobbino\\ \small Dipartimento di Matematica, Universit\`a di Pisa
\\  \small
via Filippo Buonarroti 1c, 56127 Pisa, Italy\\
\\ Margherita Solci\\ \small 
DADU, Universit\`a di Sassari\\ \small 
 piazza Duomo 6, 07041 Alghero (SS), Italy}
\date{}                                           

\begin{document}
\maketitle

\abstract \noindent We prove that a general condition introduced by Colombo and Gobbino
to study limits of curves of maximal slope allows also to characterize minimizing movements
along a sequence of functionals as curves of maximal slope of a limit functional. 

\section{Introduction}
\label{intro}
Following a vast earlier literature, the notion of minimizing movement has been introduced by De Giorgi
to give a general framework for Euler schemes in order to define a
gradient-flow type motion also for a non-differentiable energy $\phi$. 
It consists in introducing a time scale $\tau$, define a time-discrete
motion by a iterative minimization procedure in which the distance 
from the previous step is penalized in a way depending on $\tau$, 
and then obtain a time-continuous limit as $\tau\to 0$. 
This notion is at the base of modern definitions of variational motion
and has been successfully used to construct a theory of gradient flows
in metric spaces by Ambrosio, Gigli and Savar\'e \cite{AGS}. 
In particular, under suitable assumptions it can be shown that a minimizing 
movement is a curve of maximal slope for $\phi$.

When a sequence of energies $\phi_\e$ parameterized by a (small positive) 
parameter $\e$ has to be taken into account, in order to define an effective motion
one may examine the minimizing movements 
of $\phi_\e$ and take their limit  as $\e\to 0$, or, depending on the problem at hand, instead
compute the minimizing movement of the ($\Gamma$-)limit $\phi$ of $\phi_\e$. 
In general these two motions are different. 
This is due to the trivial fact that
the energy landscape of $\phi$ may not carry enough information to
describe the energy landscapes of $\phi_\e$, since local minimizers may appear or 
disappear in the limit process, as easy examples show \cite{LNM2014}.

A general approach is to proceed in the minimizing-movement scheme
letting the parameter 
$\e$ and the time scale $\tau$ tend to $0$ together. This gives a notion of
{\em minimizing movement along $\phi_\e$ at a given scale $\tau$}. In this way we can 
detect fine phenomena due to the presence of local minima. 
The limit of the minimizing movements and the minimizing movement of the limit
are recovered as extreme cases.  A first example of this approach has been given in \cite{BGN} 
for spin energies converging to a crystalline perimeter, in which case the extreme
behaviours are complete pinning and flat flow \cite{AT}. For a general choice of the parameters the
limit motion is neither of the two but depends on the ratio between $\e$ and $\tau$
and is a degenerate motion by crystalline curvature with pinning only of large sets. 
More examples can be found in \cite{LNM2014}. Note that in general the functions that 
we obtain as minimizing movements along a sequence cannot be easily rewritten as 
minimizing movements of a single functional. 

In another direction, conditions have been exhibited that ensure that the limit
of gradient flows for a family $\phi_\e$, or of curves of maximal slope, be the
gradient flow, or a curve of maximal slope, for their $\Gamma$-limit $\phi$ 
(see \cite{SS} and \cite{CG}).
This suggest that under such conditions all minimizing movements along
$\phi_\e$, at whatever scale, may converge to minimizing movements of $\phi$. 
In this paper we prove a result in that direction, showing that if a lower-semicontinuity
inequality holds for the descending slope of $\phi_\e$ then any minimizing
movement is a curve of maximal slope for $\phi$ (Theorem \ref{teo}). 
This property holds in particular in the case of convex energies (see \cite{CG}, and also \cite{AG} and \cite{LNM2014}).
Note that the limit curve of maximal slope may still depend on the
way $\e$ and $\tau$ tend to $0$, and that in the extreme case of $\e$ tending to $0$ fast enough
with respect to $\tau$ it is also a minimizing movement for the limit (Theorem \ref{teo-gamma}(b)).
In the case that all curves of maximal slope are minimizing 
movements for $\phi$, this is a kind of `commutativity result' between minimizing movements 
and $\Gamma$-convergence.   This again holds if $\phi$ is convex, which
is automatic if also $\phi_\e$ are convex. Nevertheless, in some cases it has been possible to directly prove that minimizing  movements along a sequence are minimizing movements for the 
limit, also for some non-convex energies as scaled Lennard-Jones ones
\cite{BDV}. 

The result presented in this note is suggested by the analog result for curves of maximal slope in \cite{CG},
it makes the analysis in \cite{LNM2014} more precise, and we think it will
be a useful reference for future applications. Its proof follows modifying the arguments of \cite{AGS},
which show that the minimizing movements for a single functional are curves of maximal slope, 
and is briefly presented in Section \ref{profa}.

\section{The limit result}
\label{comm}

In what follows $(X,d)$ is a complete metric space. 

\smallskip

\begin{definition}[Minimizing movements along $\phi_\e$ at scale $\tau_\e$] \label{movmin}
For all $\e>0$ let $\phi_\e\colon X\to(-\infty,+\infty]$, and let  $u_{\e}^0\in X$. 
Suppose that there exists some $\tau^*>0$ such that for every $\e>0$ and $\tau\in(0,\tau^\ast)$, 
there exists a sequence $\{u_{\e,\tau}^i\}$ which satisfies
$u_{\e,\tau}^{0}=u_{\e}^{0}$ and 
$u_{\e,\tau}^{i+1}$ is a solution of the minimum problem 
\begin{equation}
\min\Bigl\{\phi_\e(v)+\frac{1}{2\tau} d^2(v,u_{\e,\tau}^i) : v\in X\Bigr\}.
\end{equation}  
Let $\tau=\tau_\e$ be a family of positive numbers such that $\lim_{\e\to 0}\tau_\e=0$ and
define the piecewise-constant functions $\overline u_{\e}=\overline u_{\e,\tau_\e}\colon[0,+\infty)\to X$
as $\overline u_\e(t)=u_{\e,\tau}^{i+1}$ for $t\in (i\tau, (i+1)\tau]$. 
A {\em minimizing movement along $\phi_\e$ at scale $\tau_\e$} with initial data $u_{\e}^0$   
is any pointwise limit of a subsequence of the family $\overline u_\e$.
\end{definition}

From now on, we will make the following hypotheses, which ensure 
the existence of minimizing movements along the sequence $\phi_\e$ at any given scale $\tau_\e$
\cite{LNM2014}: ($u^\ast\in X$ is an arbitrary given point)

{\rm (i)} for all $\e>0$ $\phi_\e$ is lower semicontinuous; 

{\rm (ii)} there exist $C^\ast>0$ and $\tau^\ast>0$ such that
$\inf\bigl\{\phi_\e(v)+\frac{1}{2\tau^\ast}d(v,u^\ast): v\in X\bigr\}\geq C^\ast>-\infty$ for all $\e>0$;

{\rm (iii)} for all $C>0$ there exists a compact set $K$ such that   
$\{u: d^2(u,u^\ast)\leq C, |\phi_\e(u)|\leq C\}\subset K$ 
for all $\e>0$.

\begin{remark}\label{rem0} 
\rm (a) Alternatively, in the definition above we can suppose $\tau\to 0$ and choose $\e=\e_\tau\to 0$. In this way we define {\em minimizing movements along $\phi_{\e_\tau}$ at scale $\tau$};
  
(b) 
If $\phi_\e=\phi$ for all $\e$, then a minimizing movement along $\phi_\e$ at any scale is a (generalized) {\em minimizing movement for $\phi$} as defined in \cite{AGS}. \end{remark} 

For $\phi\colon X\to (-\infty,+\infty]$ we define the {\em (descending) slope} of $\phi$ as
\begin{equation}
|\partial \phi|(x)=\begin{cases}
\displaystyle\limsup_{y\to x}\frac{(\phi(x)-\phi(y))^+}{d(x,y)}&\hbox{ if $\phi(x)<+\infty$ and $x$ is not isolated}\cr
0&\hbox{ if $\phi(x)<+\infty$ and $x$ is isolated}\cr
+\infty&\hbox{ if  $\phi(x)=+\infty$.}\end{cases}\end{equation}
We say that $v\colon[0,T]\to X$ belongs to $AC^2([0,T];X)$ if there exists $A\in L^2(0,T)$ such that 
\begin{equation}
d(v(s), v(t))\leq \int_s^t A(r)\, dr \quad \hbox{ for any } 0\leq s\leq t\leq T. 
\end{equation}
The smallest such $A$ is the {\em metric derivative} of $v$ and it is denoted by $|v^\prime|$. 

\begin{definition}[Curve of maximal slope] 
A {\em curve of maximal slope} for $\phi\colon X\to (-\infty,+\infty]$ in $[0,T]$ is 
$u\in AC^2([0,T];X)$ such that there exist a non-increasing 
function $\varphi\colon[0,T]\to\mathbb R$ such that 
 \begin{equation}
\varphi(s)-\varphi(t)\geq \frac{1}{2}\int_s^t |u^\prime|^2(r)\, dr + \frac{1}{2} \int_s^t |\partial\phi|^2(u(r))\, dr 
\quad \hbox{ for any } \quad 0\leq s\leq t\leq T 
 \end{equation}
 and a Lebesgue-negligible set $E$  such that $\phi(u(t))=\varphi(t)$ in $[0,T]\setminus E$.   
\end{definition}

\begin{theorem}[Minimizing movements and curves of maximal slope]\label{teo} 
Let $\phi_\e,\phi\colon X\to(-\infty,+\infty]$ satisfy {\rm(i)-(iii)} and 
the following condition: 

{\rm (H)} for all subsequences $\phi_{\e_n}$ and $v_n\to v$ with  $\sup_n\{|\phi_{\e_n}(v_n)|
+|\partial\phi_{\e_n}|(v_n)\}<+\infty$, we have 
\begin{equation}
\label{CG} \lim_{n\to+\infty}\phi_{\e_n}(v_n)=\phi(v) \quad \hbox{ and } \quad 
 \liminf_{n\to+\infty}|\partial\phi_{\e_n}|(v_n)\geq |\partial\phi|(v).
\end{equation}   
Let $u_{\e,\tau}^0$ be such that there exist $S^\prime, S$ for which   
$d^2(u_{\e,\tau}^0,u^\ast)\leq S^\prime<+\infty$ and $|\phi_\e(u_{\e,\tau}^0)|\leq S<+\infty$. 
Then for all $\tau_\e$ any minimizing movement along $\phi_\e$ at scale $\tau_\e$ 
 with initial data $u_{\e,\tau_\e}^0$
is a curve of maximal slope for $\phi$. 
\end{theorem}

\begin{remark}\label{rem} 
{\rm (a) Note that we do not suppose that $\phi_\e$ $\Gamma$-converge to $\phi$. This is usually deduced at points with finite slope for $\phi$;

(b) The results can be proved under the more general assumption that 
hypotheses (i) and (iii) hold with respect to a weaker topology 
compatible with the distance $d$ (see \cite[Sect.~2.1]{AGS});

(c) If the initial slopes satisfy the additional bounds $|\partial \phi_\e|(u_{\e,\tau}^0)\leq S<+\infty$, then any limit curve $u$ as in the statement of Theorem~\ref{teo} satisfies
$\phi(u(t)) \leq \phi(u(0))$  for almost all $0\leq t\leq T$.}
\end{remark} 

\begin{corollary}[Curves of maximal slope and $\Gamma$-convergence]\label{cocco}
Let  $\phi_\e$ satisfy {\rm(i)-(iii)} and 
\begin{equation}\label{cono}
\phi_\e(y)\geq \phi_\e(x)-d(x,y)|\partial\phi_\e|(x) \quad \hbox{ for any } y\in X 
\end{equation}
for any $x$ such that $\phi_\e(x)$ and $|\partial\phi_\e|(x)$ are finite 
{\em (Slope Cone Property \cite{CG})}. Let $\phi_\e$ $\Gamma$-converge to $\phi_0$ 
with respect to $d$. Then the claim of Theorem {\rm\ref{teo}} holds with $\phi=\phi_0$. 
\end{corollary}
The proof of the corollary follows from \cite[Proposition 3.4]{CG}. 
Note that Corollary \ref{cocco} holds even if assumption (\ref{cono}) is weakened by subtracting in the right-hand side a lower-order term with respect to $d(x,y)$, as described in \cite[Remark~3.5]{CG}.
\begin{remark}[Convex energies] 
{\rm If $\phi_\e$ are convex, then they satisfy condition (\ref{cono}).}
\end{remark}

We note that in general the claim of Theorem \ref{teo} does not hold under the only hypothesis of $\Gamma$-convergence of $\phi_\e$. More precisely, the following theorem gives a connection between minimizing movements along a sequence 
and the (generalized) minimizing movements of the limits.  

\goodbreak
\begin{theorem}[$\Gamma$-convergence and minimizing movements]\label{teo-gamma}
Let $\phi_\e$ satisfy {\rm (i)-(iii)}. Then

{\rm(a)} there exists $\overline \tau=\overline \tau_\e$ such that 
if $\tau_\e\leq\overline \tau_\e$ then each minimizing movement along $\phi_\e$ at scale $\tau_\e$ 
is (up to subsequences) the limit as $\e\to 0$ of the (generalized) minimizing movements for $\phi_\e$;  

{\rm(b)} if $\phi_\e$ $\Gamma$-converge to $\phi$ with respect to $d$, then there exists $\overline \e=\overline \e_\tau$ 
such that if $\e_\tau\leq\overline \e_\tau$ then each minimizing movement 
along $\phi_{\e_\tau}$ at scale $\tau$ 
is (up to subsequences) a (generalized) minimizing movement for $\phi$. 
\end{theorem}
In \cite[Theorem 8.1]{LNM2014} it is proved the existence of families $\overline \tau_\e$ and $\overline \e_\tau$ 
such that the conclusions of (a) and (b) hold for minimizing movement along $\phi_\e$ exactly at scale $\overline \tau_\e$ and for 
minimizing movement exactly along $\phi_{\overline \e_\tau}$ at scale $\tau$ respectively. 
Theorem \ref{teo-gamma} follows by noticing that the arguments of that proof imply 
that we may define $\overline \tau_\e$ and $\overline \e_\tau$ such that this hold for $\tau_\e\le \overline \tau_\e$ and $\e_\tau\le\overline \e_\tau$, respectively.

In case (b), note that if the limit $\phi$ satisfies some additional differentiability hypotheses 
(see \cite[Theorems 2.3.1 and 2.3.3]{AGS}) then 
each minimizing movement is a curve of maximal slope for the limit $\phi$.

\section{Proof of Theorem \ref{teo}}\label{profa}
In order to prove a priori uniform estimates we introduce some definitions in analogy to those
in \cite[Section 3.2]{AGS}. 
For any $\e>0$ and $\delta\in(0,\tau^\ast)$, given $u\in X$ we denote by $J_{\e,\delta}(u)$ 
the set of $v$ minimizing $\phi_\e(\cdot)+ \frac{1}{2\delta} d^2(\cdot,u)$.   
Then $\tilde u_{\e,\tau}\colon [0,+\infty)\to X$ is any interpolation of the values of $\overline u_{\e,\tau}$ in $\tau\mathbb N$ 
satisfying $\tilde u_{\e,\tau}(t)\in J_{\e,t-i\tau}(u_{\e,\tau}^{i})$ for $t\in (i\tau,(i+1)\tau],$ and 
\begin{equation}\label{def-G}
G_{\e,\tau}(t)=\frac{\sup \{d(v,u_{\e,\tau}^{i}) : v \in J_{\e,t-i\tau}(u_{\e,\tau}^{i})\}}{t-i\tau} \ \ \hbox{ for } \ \  t\in (i\tau, (i+1)\tau].
\end{equation}
Note that for $t>0$ 
\begin{equation}\label{G-partialphi}
G_{\e,\tau}(t)\geq |\partial\phi_\e|(\tilde u_{\e,\tau}(t)). 
\end{equation}
Furthermore, $|u_{\e,\tau}^\prime|$ denotes the piecewise-constant function defined by 
$$|u_{\e,\tau}^\prime|(t)=\tau^{-1}d(u_{\e,\tau}^{i+1},u_{\e,\tau}^i)\hbox{ for }t\in  (i\tau, (i+1)\tau).$$

The following lemma is the analog of \cite[Lemma 3.2.2]{AGS} for a sequence $\phi_\e$.  

\begin{lemma}[A priori uniform estimates] 
Let $\e,\tau>0$ and $\tau<\tau^\ast$. Let $\phi_\e$ satisfy {\rm(i)--(iii)}, and 
$u_{\e,\tau}^0$ be such that there exist $S^\prime, S$ for which   
$d^2(u_{\e,\tau}^0,u^\ast)\leq S^\prime<+\infty$ and $|\phi_\e(u_{\e,\tau}^0)|\leq S<+\infty$. 
Let $\{u_{\e,\tau}^i\}_i$ be defined as in Definition {\rm\ref{movmin}} with initial data $u_{\e,\tau}^0$. 
Then for any $i,j\in\mathbb N$ with $i<j$ 
\begin{equation}\label{uguale} 
\phi_\e(u^i_{\e,\tau})-\phi_\e(u^j_{\e,\tau})=
\frac{1}{2}\int_{i\tau}^{j\tau} |u_{\e,\tau}^\prime|^2\, dt + 
\frac{1}{2}\int_{i\tau}^{j\tau} |G_{\e,\tau}|^2\, dt. 
\end{equation}
Moreover, fixed $T>0$, there exists a constant $C$ depending only on $T,S,S^\prime,\tau^\ast,u^\ast$ such that, for 
$\tau<\tau^\ast/8$, for any 
$N$ with $N\tau\leq T$ 
\begin{eqnarray}
&&\label{stima-dist}d^2(u^N_{\e,\tau},u^\ast)\leq C, \quad |\phi_\e(u^N_{\e,\tau})|\leq C\\
&&\label{stima-tilde}d^2(\tilde u_{\e,\tau}(t), \overline u_{\e,\tau}(t))\leq C\tau \quad \hbox{ in } \ [0,T]\\
&&\label{stima-a}\int_0^{N\tau}|u^\prime_{\e,\tau}|^2(r)\, dr\leq \phi_\e(u^0_\e)-\phi_\e(u^N_{\e,\tau})\leq C\\
&&\label{stima-b}\int_0^{N\tau}|G_{\e,\tau}|^2(r)\, dr\leq \phi_\e(u^0_\e)-\phi_\e(u^N_{\e,\tau})\leq C. 
\end{eqnarray}
\end{lemma}
The proof follows that of \cite[Lemma 3.2.2]{AGS}, where a single $\phi$ is considered. 
The uniform bounds are ensured by the condition (ii) for the sequence $\phi_\e$ and by the 
uniform bounds on the initial data, which allow to prove that for any $i$ 
\begin{equation} 
d^2(u_{\e,\tau}^i, u^\ast)
\leq  2S^\prime+2\tau^\ast S
-2\tau^\ast C^\ast+\frac{4}{\tau^\ast}\sum_{j=1}^i \tau \, d^2(u_{\e,\tau}^j,u^\ast). 
\end{equation}
An application of the Gronwall Lemma to this inequality gives the required uniform bounds 
(see \cite[pp.~67--68]{AGS}).

\smallskip

With the aid of this lemma we can prove Theorem \ref{teo} by showing that, up to subsequences, 
$v_n=\tilde u_{\e_n,\tau_{\e_n}}(t)$ satisfies the hypotheses of condition (H) for almost all $t$.   
Estimate (\ref{stima-a}) implies the weak convergence in $L^2_{\rm loc}[0,+\infty)$ 
of a subsequence $|u^\prime_{\e_n,\tau_n}|\wto A$. 
Since each $\phi_{\e_n}(\tilde u_{\e_n,\tau_n}(\cdot))$ is not increasing (see \cite[Lemma 3.1.2]{AGS}), 
we may apply Helly's Lemma (e.g. \cite[Lemma 3.3.3]{AGS}) obtaining  
(up to a further subsequence)
\begin{equation}
\varphi(t)=\lim_{n\to +\infty}\phi_{\e_n}(\tilde u_{\e_n,\tau_n}(t)) \quad \hbox{ for any } \ t\in [0,T];
\end{equation}
in particular, the family $|\phi_{\e_n}(\tilde u_{\e_n,\tau_n}(t))|$ is equibounded for any $t\in [0,T]$. 
From (\ref{G-partialphi}), estimate (\ref{stima-b}) and Fatou's Lemma 
$$\int_{0}^{T} \left( \liminf_{n\to +\infty} |\partial\phi_{\e_n}|(\tilde u_{\e_n,\tau_n}(t))\right)\,dt\leq \liminf_{n\to
+\infty}\int_{0}^{T} |\partial\phi_{\e_n}|^2(\tilde u_{\e_n,\tau_n}(t)) \,dt<+\infty,$$
we get that 
there exists a Lebesgue-negligible set $E$ such that  
\begin{equation}
	 \liminf_{n\to +\infty} |\partial\phi_{\e_n}|(\tilde u_{\e_n,\tau_n}(t))<+\infty
	\quad\quad
	\forall t\in[0,T]\setminus E.
	\label{est:slliminf}
\end{equation}
The estimates in (\ref{stima-dist}) and the compactness hypothesis (iii)  
imply that the family $\overline u_{\e_n,\tau_n}$ is contained in a compact subset of $X$. 
We set $s(n)= i_n\tau_n$ with $s(n)\leq s$ and $s(n)\to s$, and  
$t(n)= j_n\tau_n$ with $t(n)\geq t$ and $t(n)\to t$. 
By the definition of $|u_{\e_n,\tau_n}^\prime|$ and by weak convergence we get 
\begin{equation}\label{stima-A}
\limsup_{n\to +\infty}d(\overline u_{\e_n,\tau_n}(s), \overline u_{\e_n,\tau_n}(t)) 
 \leq \limsup_{n\to +\infty}\int_{s(n)}^{t(n)}|u_{\e_n,\tau_n}^\prime|(r)\, dr
\leq\int_{s}^{t}A(r)\, dr
\end{equation}
for any $s,t\in [0,T]$. This estimate and the pre-compactness of  
the family $\overline u_{\e_n,\tau_n}$ allow to use a Ascoli-Arzel\`a argument (see e.g. \cite[Proposition 3.3.1]{AGS})  
to obtain the existence of a function $u$ pointwise limit of a further subsequence of 
$\overline u_{\e_n,\tau_n}(t)$. 
Since (\ref{stima-tilde}) holds, also $\tilde u_{\e_n,\tau_n}(t)$ converges to $u(t)$ for any $t\in [0,T]$.
Since (\ref{est:slliminf}) is in force, we can extract for any $t \in [0,T] \setminus E$
 a ($t$-dependent)
sequence $n_{k}\to +\infty$ such that
$$
	\lim_{k\to +\infty}|\partial\phi_{\e_{n_k}}|(\tilde u_{\e_{n_k},\tau_{n_k}}(t))=
	\liminf_{n\to +\infty}|\partial\phi_{\e_n}|(\tilde u_{\e_n,\tau_n}(t))<+\infty.
$$
Applying (H) to $\tilde u_{\e_{n_k},\tau_{n_k}}(t)$, for $t\in [0,T]\setminus E$ we get 
\begin{eqnarray}
&&\label{CG1} \phi(u(t))=\lim_{k\to +\infty}\phi_{\e_{n_k},\tau_{n_k}}(\tilde u_{\e_{n_k},\tau_{n_k}}(t)) = \lim_{n\to +\infty}\phi_{\e_n,\tau_n}(\tilde u_{\e_n,\tau_n}(t)) =\varphi(t)\\
&&\label{CG2}  |\partial\phi|(u(t)) \leq \lim_{k\to +\infty}|\partial\phi_{\e_{n_k}}|(\tilde u_{\e_{n_k},\tau_{n_k}}(t)) = \liminf_{n\to +\infty} |\partial\phi_{\e_n}|(\tilde u_{\e_n,\tau_n}(t)).
\end{eqnarray}
It remains to prove that $u$ is a curve of maximal slope for $\phi$. 
Passing to the limit in (\ref{stima-A}) 
it follows that $u\in AC^2([0,T];X)$ and that for almost all $t\in[0,T]$ 
$|u^\prime|(t)\leq A(t)$. The weak convergence ensures 
that 
\begin{equation}\label{uprime}\liminf_{n\to+\infty}\int_{s(n)}^{t(n)} |u_{\e_n,\tau_n}^\prime|^2(r)\, dr
\geq \int_{s}^{t} A^2(r)\, dr\geq  \int_{s}^{t} |u^\prime|^2(r)\, dr. 
\end{equation}
Estimates (\ref{G-partialphi}) and (\ref{CG2}), and an application of Fatou's Lemma imply 
\begin{equation}\label{G}
\displaystyle\liminf_{n\to+\infty}\int_{s(n)}^{t(n)} |G_{\e_n,\tau_n}|^2(r)\, dr
\geq  \displaystyle\liminf_{n\to+\infty}\int_{s(n)}^{t(n)} |\partial\phi_{\e_n}|^2(\tilde u_{\e_n,\tau_n}(r))\, dr
\geq  \displaystyle \int_{s}^{t} |\partial\phi|^2(u(r))\, dr. 
\end{equation}
By the non-increasing monotonicity of $\phi_{\e_n}(\tilde u_{\e_n,\tau_n}(\cdot))$ we have the inequalities 
$\phi_{\e_n}(u^{i_n}_{\e_n,\tau_n})\geq \phi_{\e_n}(\tilde u_{\e_n,\tau_n}(s))$ and 
$\phi_{\e_n}(u^{j_n}_{\e_n,\tau_n})\leq \phi_{\e_n}(\tilde u_{\e_n,\tau_n}(t))$, 
so that, also using (\ref{uprime}), (\ref{G}) and (\ref{CG1}), 
from (\ref{uguale}) 
we obtain 
\begin{equation}\label{tesi} 
\varphi(s)-\varphi(t)\geq \frac{1}{2}\int_s^t |u^\prime|^2 \, dr 
+ \frac{1}{2} \int_s^t |\partial \phi|^2(u(r)) \, dr, 
\end{equation} 
where $\phi(u(t))=\varphi(t)$ in $[0,T]\setminus E$.  

\section*{Acknowledgements}
The authors thank the anonymous referee for her/his useful comments. The second author has been partially supported by Dr. Max R\"ossler, the Walter Haefner Foundation and the ETH Zurich Foundation.

\end{document}